\documentclass[]{amsart}
\usepackage{amssymb,amsxtra}
\theoremstyle{plain}
\newtheorem{thm}{Theorem}[section]

\newtheorem{cor}[thm]{Corollary}
\newtheorem{defn}[thm]{Definition}
\newtheorem{prop}[thm]{Proposition}
\theoremstyle{remark}
\newtheorem{rem}[thm]{Remark}
\newtheorem{ex}[thm]{Example}
\newtheorem{exs}[thm]{Examples}
\title[KMS states]{KMS states on C*-algebras associated to expansive maps}
\author[Kumjian]{Alex Kumjian}
\address{Department of Mathematics, University of Nevada, Reno NV
89557, USA}
\email{alex@unr.edu}
\author[Renault]{Jean Renault}
\address{D\'epartment de Math\'ematiques, Universit\'e d'Orl\'eans,
45067 Orl\'eans, France}
\email{renault@labomath.univ-orleans.fr}
\date{27 Oct 2003}
\begin{document}
\begin{abstract} 
Using Walters' version of the Ruelle-Perron-Frobenius Theorem we show the
existence and uniqueness of KMS states for a certain one-parameter
group of automorphisms on a C*-algebra associated to a positively
expansive map on a compact metric space.
\end{abstract}
\maketitle

\section{Introduction}\label{intro}
In the C*-algebraic formulation of quantum statistical mechanics KMS states
are equilibrium states associated to the one-parameter group of
automorphisms of time evolution (see \cite[\S 5.3]{br}).  This
association mirrors that of a 
faithful normal state on a von Neumann algebra to its modular
automorphism group.  Our main result concerns the C*-algebra associated
to a positively expansive exact local homeomorphism $T$ defined on a compact
metric space $X$ (see \cite{d,r3}).  A continuous real-valued function
$\varphi$ on $X$ gives rise to a one-parameter group of
automorphisms on the C*-algebra.  We apply Walters' version of the
Ruelle-Perron-Frobenius Theorem and groupoid techniques to 
establish existence and uniqueness of KMS states under mild hypotheses 
on $\varphi$; in particular, if $\varphi$ is positive and H\"older
continuous there is a unique KMS state (this generalizes recent work
of Exel, cf.\ \cite{ex3}). 

When a C$^*$-algebra $A$ is written as the C$^*$-algebra of a groupoid $G$
(by the
construction of \cite{r1}),  a natural problem is to describe its structure
in terms of
$G$. For example, positive results about the ideal structure (see
\cite{r2}) or the
K-theory (see
\cite{tu}) of groupoid C$^*$-algebras have been obtained, especially in the
amenable case.
The scope of this paper is more limited. Let $G$ be a locally compact
groupoid with Haar
system (for convenience, we assume that $G$ is \'etale and Hausdorff) and
let $C^*(G)$ be
its C$^*$-algebra. Every continuous cocycle $c\in Z^1(G,{\bf R})$ defines a
one-parameter automorphism group $\alpha_c$ of $C^*(G)$ as in
\cite[Section II.5]{r1}. First, we show in
Proposition
\ref{kms} that if $c^{-1}(0)$ is principal, then the
KMS states of $\alpha_c$ at inverse temperature $\beta$ are exactly the
quasi-invariant
probability measures which admit $e^{-\beta c}$ as Radon-Nikodym derivative.
This extends
\cite[Proposition II.5.4]{r1}, which gave a complete answer only when $G$
itself is
principal. There certainly exist improvements of this result (see Remark
\ref{no atoms}),
but it suffices for our purpose. Then, we apply this result to the
Cuntz-like
C$^*$-algebras of
\cite{d,r3}. These algebras are of the form $C^*(G)$, where $G=G(X,T)$ and
$(X,T)$
is a dynamical system where
$X$ is a compact space and $T$ is a surjective local homeomorphism (see
below).
Proposition \ref{kms} reduces the KMS problem to a Radon-Nikodym problem
which fits the
framework of the thermodynamical formalism. Under classical assumptions on
the dynamical
system $(X,T)$ and on the potential $\varphi\in C(X,{\bf R})$ defining the
cocycle $c_\varphi\in Z^1(G,{\bf R})$, Ruelle's Perron-Frobenius
theorem can be applied to obtain existence and
uniqueness results. Theorem \ref{main} summarizes our results about KMS
states of these
automorphism groups. We assume that $T$ is positively expansive and exact
(see Proposition \ref{ws=min}). We first give a necessary and
sufficient condition for the existence of KMS states at inverse temperature
$\beta$,
namely the equation
$P(T,-\beta\varphi)=0$, where $P$ is the pressure. Under the assumption that
$\varphi$
satisfies the Walters condition, we give a necessary condition 
($c_\varphi^{-1}(0)$ is principal)
for the existence of KMS states. If $\varphi$
satisfies the Bowen condition and has constant sign, then we have
the existence
and the uniqueness of the KMS inverse temperature and state. Our results,
which complete
\cite[Section 4.2]{r4}, strictly cover previous results about KMS states of
gauge
one-parameter automorphism groups of Cuntz and Cuntz-Krieger algebras
obtained by a number
of authors (e.g. \cite{op, ev, efw, z, kp, ex2}); although they do not use
groupoid
techniques as we do, Kerr-Pinzari in \cite{kp} and Exel in \cite{ex2}
also make essential use of Ruelle's transfer operator.

Our work is organized as follows. In Section 2, we recall the definition of
the groupoid
$G(X,T)$ and some facts pertaining to the dynamics of a positively expansive
map on a
compact metrice space. Section 3 contains our results about KMS states and a
list of
examples. Besides the cases previously studied, we give an example of a
unique KMS state
for a potential which does not have a constant sign. Section 4 gives some
properties of
the associated crossed product C$^*$-algebras.

The first author wishes to thank the Universit\'e d'Orl\'eans for its
support.

\section{Preliminaries}\label{pre}
For $X$ a locally compact space and $T : X \rightarrow X$ a local
homeomorphism,  we let $G = G(X, T)$ denote the associated groupoid
(see \cite{d,r3}) where
$$
G(X, T) = \{ (x, m - n, y) : T^m x = T^n y \}.
$$ 
This is an \'etale groupoid (a topological groupoid with structure
maps that are local homeomorphisms); note that the unit space $G^0$
may be identified with $X$ via the map $x \mapsto (x, 0, x)$.   The
subgroupoid for which  $m = n$, denoted $R = R(X, T)$, is an
approximately proper equivalence relation (see \cite{r4,r5}).
Note that $R = \cup_n R_n$ where $(x, y) \in  R_n$ if
$T^m x = T^n y$ and $C^*(R_n)$ is strongly Morita equivalent to
$C_0(X)$.

Let $X$ be a compact metric space with metric $\rho$ which is fixed
throughout and assumed to contain infinitely many points.  A
continuous surjection $T : X \rightarrow X$ is said to be
{\em positively expansive} if there is an $\epsilon > 0$ such that if
$x \ne y$ then $\rho(T^n x, T^ny) \ge \epsilon$ for some $n \ge 0$.
Note that this property does not depend on the choice of the metric and
can be expressed topologically as follows: there is an open subset $U$
of $X \times X$ which contains the diagonal such that for all $x \not= y$
there is an $n \in {\bf N}$ with $(T^nx, T^ny) \notin U$.
Observe that $T$ is not a homeomorphism, for we have assumed that
$X$ is not finite (see \cite[2.2.12]{ah}).

Let $T : X \rightarrow X$ be positively expansive.
By \cite{rd} there are an equivalent metric $\rho'$ and constants
$\tau > 0$, $\kappa > 1$ such that if $\rho'(x, y) < \tau$ then
$\rho(T x, Ty) \ge \kappa\rho'(x, y)$.
If $T$ is open then it is a local homeomorphism (such a map is called {\em
expanding} in \cite{ah}); see \cite[2.2.21]{ah} for an example due to
Rosenholtz of a positively expansive map which is not open and
therefore not a local homeomorphism.

Then
Let $T : X \rightarrow X$ be a local homeomorphism.  Then $T$ is said to be
{\em exact} if for every non-empty open set $U \subset X$ there is an
$n > 0$ such that $T^n(U) = X$ (this property is called mixing in  \cite[\S
2]{fj}).
We say that $T$ satisfies the {\em weak specification condition}
if for every $\epsilon > 0$ there is $n > 0$ such that for
all $x \in X$, $T^{-n} x$ is $\epsilon$-dense, that is, we have
$$
X = \bigcup_{y \in T^{-n} x} B(y, \epsilon),
$$
where $T^{-n} x = \{ y \in X \mid T^n y = x \}$ and $B(y, \epsilon)$ is the
open ball of radius $\epsilon$ centered at $y$ (cf.\ \cite[1.2]{w3} where
this is listed as one of four equivalent conditions for a positively
expansive local homeomorphism).

\begin{prop}\label{ws=min}
Let $T : X \rightarrow X$ be a local homeomorphism.  Then the following
conditions are equivalent:
\begin{itemize}
\item[(i)] $R$ is minimal.
\item[(ii)] $T$ is exact.
\item[(iii)] $T$  satisfies the weak specification condition.
\end{itemize}
\end{prop}
\begin{proof}
The equivalence of (i) and (ii) follows by the compactness of $X$.
For the minimality of $R$ is equivalent to the condition that the
saturation of any non-empty open set is all of $X$.  Now suppose that
(ii) holds and fix $\epsilon > 0$.  By compactness there are
$x_1, \dots, x_m \in X$ such that $X \subset \cup_i B(x_i, \epsilon/2)$.
By (ii) there is $n > 0$ such that $T^n(B(x_i, \epsilon/2)) = X$ for
$i = 1, \dots, m$.  To show that $T$ satisfies the weak specification
condition, we verify that $T^{-n} x$ is $\epsilon$-dense for every
$x \in X$.  Given $x, z \in X$ there is an $i$ such that
$\rho(z, x_i) < \epsilon/2$.  Since $T^n(B(x_i, \epsilon/2)) = X$,
$B(x_i, \epsilon/2) \cap T^{-n} x \not= \emptyset$.
Choose an element $y$ in this set; then we have
$$
\rho(z, y) \le \rho(z, x_i) + \rho(x_i, y) < \epsilon
$$
and $T^{-n} x$ is $\epsilon$-dense as required.
Now suppose that $T$  satisfies the weak specification condition.  We
verify that $R$ is minimal, that is, every orbit under $R$ is dense.
Given $x \in X$, the orbit of $x$ is given by $\cup_nT^{-n}(T^n x)$;
for $\epsilon > 0$, there is $n > 0$ so that $T^{-n}(T^n x)$ is
$\epsilon$-dense 
by weak specification.  Hence, the orbit of $x$ is dense.
\end{proof}
\begin{ex}\label{ber}
Fix $n \in {\bf N}$, $n > 1$, and let
$X = \{(x_0, x_1, \dots) \mid x_i = 1, 2, \dots, n  \}$.
For $x, y \in X$ with $x \not= y$, set $\rho(x, y) = 2^{-j}$ where
$j = \min\{i \mid x_i \not= y_i \}$ and define $T : X \rightarrow X$ by
$T(x_0, x_1, \dots) = (x_1, x_2, \dots)$.  Then $X$ is a compact
metric space and $T$ is a positively expansive local homeomorphism
which satisfies the weak specification condition.  The map $T$ is
often called the Bernoulli shift on $n$ letters.  Note that
$C^*(G)$ is the Cuntz algebra $\mathcal{O}_{n}$ (see \cite{cu}),
which is generated by $n$ isometries $S_1, \dots, S_n$ subject
to the relation 
$$
1 = \sum_{j=1}^n S_jS_j^*.
$$
\end{ex}
\begin{ex}\label{sft}(see \cite{ck})
Let $A$ be an irreducible $n \times n$ zero-one matrix
$A = (A(i,j))_{1 \le i,j \le n}$.  Then $A$ defines a closed subset
$X_A$ of $X$ above 
$$
X_A = \{ x \in X : A(x_i, x_{i+1}) = 1 \}.
$$
Note that $TX_A \subset X_A$ and let $T_A: X_A \rightarrow X_A$ denote the
restriction.  Such a map is called a subshift of finite type.
Note that $T_A$ is a positively expansive local homeomorphism;
but it does not necessarily  satisfy the weak specification condition
unless $A$ is primitive (i.e.\ there is a positive integer $k$ such
that all entries of $A^k$ are nonzero).
Note, $C^*(G) \cong {\mathcal O}_A$, where ${\mathcal O}_A$ is the
Cuntz-Krieger algebra
associated to the matrix $A$; ${\mathcal O}_A$  is also generated by $n$
isometries. 
\end{ex}
\begin{ex}\label{tor}
Fix $n > 1$; let $X = {\bf T} = {\bf R}/{\bf Z}$ and define $T: X
\rightarrow X$ by 
$Tx = nx$.  We endow $X$ with the arclength metric.  Then $T$ is a
positively expansive local homeomorphism which satisfies the weak
specification condition.   More generally an integer $k \times k$
matrix, all of whose eigenvalues exceed one in absolute value, defines
a positively expansive local homeomorphism on $X = {\bf T}^k = {\bf
R}^k/{\bf Z}^k$
(see \cite{ah}). 
\end{ex}

Let $C(X ; {\bf R})$ denote the space of continuous real-valued functions
on ${\bf R}$.  Then for $\varphi \in C(X ; {\bf R})$ we define the transfer
operator
$\mathcal{L}_\varphi : C(X ; {\bf R}) \rightarrow C(X ; {\bf R})$ by the
formula
$$
(\mathcal{L}_\varphi f)(x) = \sum_{y \in T^{-1} x} e^{\varphi(y)}f(y).
$$
The following proposition is an easy consequence of Theorem 1.3 and
Corollary 2.3 of \cite{w3}. The pressure $P(T, \varphi)$ of $T$ at $\varphi$
is defined in \cite[\S 9.1]{w2}.

\begin{prop}\label{press}
Let $T : X \rightarrow X$ be a local homeomorphism which is positively
expansive and exact
and let $\varphi \in C(X ; {\bf R})$.   There is a unique 
$\lambda > 0$ such that
$$
\mathcal{L}_\varphi^*\mu = \lambda\mu
$$
for some probability measure $\mu$ on $X$.
Moreover, $\log \lambda = P(T, \varphi)$ and the sequence
$\frac{1}{n}\log(\mathcal{L}_\varphi^n1)$
converges uniformly to the constant $P(T, \varphi)$.
\end{prop}
Note that $P(T, 0) = h(T)$, the topological entropy of $T$ (see
\cite[9.7i]{w2}).
\begin{cor}\label{h(T)>0}
Let $T : X \rightarrow X$ be as above.  Then $h(T) > 0$.
\end{cor}
\begin{proof}
Apply the above proposition with $\varphi = 0$.  Then there is a
probability measure $\mu$ on $X$ such that
$$
\mathcal{L}_0^*\mu = e^{h(T)}\mu.
$$
Arguing as in \cite[Corollary 4.3]{r3} we may deduce that $\mu$ is
quasi-invariant under $G$.  By Proposition \ref{ws=min}, $R$ is
minimal and therefore $\mu$ must be faithful.  We have
$\mathcal{L}_0(1) > 1$.  Hence
$$
e^{h(T)} = \mathcal{L}_0^*\mu(1) = \mu(\mathcal{L}_0(1)) > \mu(1) = 1.
$$
Thus $h(T) > 0$ as required.
\end{proof}
To ensure that the measure $\mu$ is unique we will need to
impose a condition on the function $\varphi$.

\begin{defn} (see \cite[\S 1]{w3}, \cite{ru})
A function $\varphi \in C(X ; {\bf R})$ is said to satisfy the {\em
Bowen
condition} with respect to $T : X \rightarrow X$ if there are $\delta, C >
0$ such
that 
$$
\sum_{i = 0}^{n - 1} \varphi(T^ix) - \varphi(T^iy) \le C
$$
for all $x, y \in X$ and $n > 0$ such that
$d(T^ix, T^iy) \le \delta$ for $0 \le i \le n - 1$.
\end{defn}

The {\em Walters condition} introduced in \cite{w1} by P.~Walters
requires that $C$ can be 
chosen arbitrarily small (then $\delta$ depends on $C$). Recall that a real-valued
function
$\varphi$ defined on a compact metric space $(X, \rho)$ is said to be
H\"older continuous
if there are constants $r, K > 0$ such that $|\varphi(x) - \varphi(y)| \le
K\rho(x, y)^r$
for all $x, y \in X$.  Observe that a H\"older continuous function satisfies
the Walters condition. To show this we may assume
that we are dealing with a Reddy metric; if $\delta \le \tau$, then the
above sum is not
greater than
$C =K\delta^r\kappa^r/(\kappa^r-1)$. If $X$ is a smooth manifold, then
differentiability implies H\"older continuity.

By Walters'
 version of the Ruelle-Perron-Frobenius Theorem (cf.\
\cite[Theorem 2.16]{w3} and \cite[Theorem 8]{w1}) we have:

\begin{thm}\label{wrpf}
Let $T : X \rightarrow X$ be a 
local homeomorphism which is positively expansive and exact 
and let $\varphi \in C(X ; {\bf
R})$.
\begin{itemize}
\item[(i)] Assume that $\varphi$
satisfies the Bowen condition.
Then there are unique $\lambda > 0$ and $\mu \in M(X)$ such that
$$
\mathcal{L}_\varphi^*\mu = \lambda\mu.
$$
\item[(ii)] Assume that $\varphi$
satisfies the Walters condition. Then there is also  $h\in C(X,{\bf R})$,
$h>0$ such that 
$$
\mathcal{L}_\varphi h= \lambda h.
$$
\end{itemize}

\end{thm}
A proof of the assertion $(i)$ of this theorem 
is given in \cite[Theorem 6.1 and Proposition 7.2]{r5}.

\section{KMS states}\label{C*}

Our main result on existence and uniqueness of KMS states is contained
in Theorem \ref{main} below.  If $G$ is an \'etale groupoid we let 
$C^*(G)$ denote the full C$^*$-algebra of $G$ as defined in 
\cite[Chapter II]{r1}.   Let $T : X \rightarrow X$ be a 
local homeomorphism which is positively expansive and exact. 
Then if  $G = G(X , T)$ is the groupoid constructed above, $C^*(G)$ is
purely infinite and simple (see \cite{claire}). We shall mainly be
concerned with KMS states on such C$^*$-algebras.

\begin{defn}
Let $A$ be a C*-algebra, $\beta \in {\bf R}$ and let $\alpha : {\bf R}
\rightarrow {\rm Aut} (A)$
be a strongly continuous action.  Then a state $\omega$ on
$A$ is said to satisfy the KMS condition for $\alpha$ at inverse
temperature $\beta$ if
$$\omega(b\alpha_{i\beta}(a))=\omega(ab)$$
for all $a,b \in A$ with $a$ entire for $\alpha$
(see \cite{bek}, \cite[\S 8.12]{p}, \cite[\S 5.3]{br}).
\end{defn}

Let $G$ be an \'etale groupoid and let $c \in Z^1(G, {\bf R})$ be a
continuous real-valued one-cocycle.  
Let $\alpha^c: {\bf R} \rightarrow {\rm Aut}(C^*(G))$ denote the
associated action on $C^*(G)$ 
defined for $f \in C_c(G)$ by
$(\alpha^c_t(f))(\gamma)=e^{itc(\gamma)}f(\gamma)$  for
$\gamma \in G$ and $t \in {\bf R}$.  Let $\mu$ be a quasi-invariant
probability
measure on the unit space $G^{(0)}$ and $\beta \in {\bf R}$; then $\mu$ is
said
to satisfy the $(c, \beta)$-KMS condition if the modular function of $\mu$
(i.e.\ the Radon-Nikodym derivative $dr^*\mu/ds^*\mu$) is $e^{-\beta c}$
(see \cite[Def.\ I.3.15]{r1}).  The associated state on $C^*(G)$
defined by $\omega_\mu(f) = \int_{G^{(0)}}f(x)\,d\mu(x)$ for
$f \in C_c(G)$ satisfies the KMS condition for $\alpha$ at inverse
temperature $\beta$. It is shown in \cite[Prop.\ II.5.4]{r1} that if $G$ is
principal then
every KMS state arises in this fashion.  Since a groupoid
associated to a local homeomorphism need not be principal, we need an
extension of this result. For the sake of simplicity, we only consider the
case of a
locally compact groupoid which is \'etale and Hausdorff.

\begin{prop}\label{kms}
Let $G$ be an \'etale and Hausdorff groupoid. Let $c \in Z^1(G, {\bf R})$ be
a continuous
real-valued one-cocycle and let $\alpha$ denote the associated
action. Assume that $c^{-1}(0)$ is principal. Then every KMS state for
$\alpha$ at inverse
temperature
$\beta$ is of the form
$\omega_\mu$ for some quasi-invariant probability measure $\mu$ on $G^{(0)}$
with
Radon-Nikodym derivative $dr^*\mu/ds^*\mu = e^{-\beta c_\varphi}$.
\end{prop}

\begin{proof}  Our assumption means that $c^{-1}(0)\cap G'\subset
X$, where
$G'$ is the isotropy bundle
$\{\gamma\in G:r(\gamma)=s(\gamma)\}$ and
 $X=G^{(0)}$. Replacing $c$ by $\beta c$, we may assume that $\beta=1$.
Thus, let $\omega$ be a KMS state for $\alpha$ at inverse temperature $1$.
Then 
the restriction of $\omega$ to $C_0(X)$ defines a probability
measure $\mu$ on $X$.  It is shown in \cite[II.5.4]{r1} that $\mu$ is
quasi-invariant with
the requisite modular function.

The restriction of $\omega$ to $C_c(G)$ is a complex Radon measure and
defines a complex
measure $\nu$ on $G$ ($\mu$ is the restriction
of $\nu$ to $X$). We show just as in \cite[II.5.4]{r1} that the support of
$\nu$ is
contained in the closed set
$G'$. Because of the KMS condition, we have $\omega((h\circ
r)f)=\omega(f(h\circ s))$ for all $f\in C_c(G), h\in C_c(X)$. Suppose that
the support of
$g\in C_c(G)$ does not meet $G'$. Then there exist open sets $U_i,V_i\subset
X$, $i\in I$ finite, such that $\{G_{U_i}^{V_i},\, i\in I\}$ covers the
support of
$g$ and $U_i\cap V_i=\emptyset$ for all $i\in I$. Using a partition of
unity, we can write
$g=\sum_I g_i$ where the support of $g_i$ contained in $G_{U_i}^{V_i}$. For
all $i\in I$,
there exist
$h_i\in C_c(X)$ such that $(h\circ r)g_i=g_i$ and $g_i(h\circ s)=0$. This
implies
$\omega(g_i)=0$ and $\omega(g)=0$.
Next we note that $\omega$ is invariant under $\alpha$. This implies that
for all $f\in
C_c(G)$ and all $t\in {\bf R}$, we have
$\omega(e^{itc}f)=\omega(f)$. This implies that for all $t\in{\bf R}$,
$e^{itc}\equiv 1$
on the support of
$\nu$. Therefore the support of $\nu$ is contained in
$c^{-1}(0)$. Thus under our assumption that $c^{-1}(0)\cap G'$ is contained
in $
G^{(0)}$, the measures $\nu$ and $\mu$ agree. Then $\omega$ and $\omega_\mu$
agree on
$C_c(G)$, hence on $C^*(G)$.
\end{proof}

Given a local homeomorphism $T : X \rightarrow X$ and  $\varphi \in C(X ;
{\bf R})$, we
construct the groupoid $G=G(X,T)$ as above. Every $\varphi \in C(X ; {\bf
R})$
defines a continuous one-cocycle
$c = c_\varphi
\in Z^1(G, {\bf R})$ by the formula
$$
c_\varphi(x, m - n, y) =
\sum_{i = 0}^{m - 1} \varphi(T^ix) - \sum_{j = 0}^{n - 1} \varphi(T^jy).
$$
Moreover every element of $Z^1(G, {\bf R})$ is of this form (see
\cite{dkm}). In that case, the condition  $c_\varphi^{-1}(0)$ is principal
means 
that 
$$T^n(x)=x\quad \hbox{and}\quad n\ge 1\qquad\Rightarrow \qquad
\varphi_n(x)\equiv\varphi(x)+\varphi(Tx)+\ldots +
\varphi(T^{n-1}x)\not= 0.$$

It implies that $\varphi$ is not of the form $\psi\circ T-\psi$. This
condition is
trivially satisfied when
$\varphi$ is strictly positive (or strictly negative).

\begin{cor}\label{rn}
Let $T : X \rightarrow X$ be a local homeomorphism and let $\varphi \in C(X
; {\bf R})$.
Let $\alpha$ denote the action associated to the cocycle $c_\varphi
\in Z^1(G, {\bf R})$.
If $c_\varphi^{-1}(0)$ is principal, then every KMS state for $\alpha$
at inverse temperature $\beta$ is of the form
$\omega_\mu$ for some quasi-invariant probability measure $\mu$ with
Radon-Nikodym derivative $dr^*\mu/ds^*\mu = e^{-\beta c_\varphi}$.
\end{cor}

\begin{rem}\label{no atoms} In some cases, one can reach the same conclusion
without any
assumption on
$\varphi$. For example, if the restriction of the KMS state to
$C(X)$ induces a probability measure on $X$ with no atoms and if $T$ is a
local homeomorphism which is positively expansive and exact, one 
also obtains that every KMS state is of the form $\omega_\mu$
(using the fact that periodic points of a given period are also isolated).
\end{rem} 

We are now ready to state the main result concerning the existence and
uniqueness of KMS states for the action associated to the cocycle
$c_\varphi$; a key step is to show that $P(T, -\beta\varphi) = 0$ has a
unique
solution (see Corollary 4.3 of \cite{r3} and following remarks).
\begin{thm}\label{main}
Let $T : X \rightarrow X$ be a 
local homeomorphism which is positively expansive and exact.
Let $\varphi \in C(X ; {\bf R})$ and let
$\alpha$ be the automorphism group associated to the cocycle
$c_\varphi \in Z^1(G, {\bf R})$.
\begin{itemize}
\item[(i)] There is a KMS state for $\alpha$ at inverse temperature
$\beta\in{\bf R}$ if and only if
$P(T,-\beta\varphi)=0$.
\item[(ii)]  Let
$\beta\in{\bf R}$ such that $P(T,-\beta\varphi)=0$. If $c_{\varphi}^{-1}(0)$
is principal
and $\varphi$ satisfies the Bowen condition, then the KMS state at
inverse temperature
$\beta$ is unique.
\item[(iii)] Assume that $c_\varphi^{-1}(0)$ is not principal and that
$\varphi$ satisfies
the Walters condition. Then, for all $\beta\in{\bf R}$,
$P(T,-\beta\varphi)>0$. In
particular,
$\alpha$ has no KMS state.
\item[(iv)] If  $\varphi$ is strictly positive (resp.\ strictly negative),
there is a unique
$\beta$ such that $P(T,-\beta\varphi)=0$ and therefore a unique  inverse
temperature $\beta$ at which KMS states exist; moreover,
$\beta > 0$ (resp.\ $\beta < 0$): indeed we have 
$$ 
\frac{h(T)}{\sup\varphi} \le \beta \le \frac{h(T)}{\inf\varphi}.
$$
If furthermore $\varphi$ satisfies the Bowen condition, the
KMS state at inverse temperature $\beta$ is unique.
\end{itemize} 
\end{thm}
\begin{proof} For $(i)$, we note that if $\omega$ is a KMS state for
$\alpha$ at inverse
temperature $\beta$, the measure $\mu$ on $X$ given by its restriction to
$C(X)$ is
quasi-invariant with Radon-Nikodym derivative $e^{-\beta c_{\varphi}}$.
This implies
(\cite[Proposition 4.2]{r3}) that
$\mathcal{L}_{-\beta \varphi}^*\mu = \mu$ and therefore, according to
Proposition \ref{press}, that $P(T,-\beta\varphi)=0$. Conversely, if
$P(T,-\beta\varphi)=0$, the same proposition gives the existence of a
quasi-invariant
probability measure $\mu$ with Radon-Nikodym derivative $e^{-\beta
c_{\varphi}}$. Then,
$\omega_\mu$ is a KMS state for $\alpha$ at inverse
temperature $\beta$. For $(ii)$, we note that the assumption that
$c_{\varphi}^{-1}(0)$
is principal guarantees that every KMS state is of the form $\omega_\mu$,
where $\mu$ is
a Perron-Frobenius eigenvector $\mathcal{L}_{-\beta \varphi}^*\mu = \mu$ and
the second
assumption guarantees the uniqueness of this Perron-Frobenius eigenvector.
For $(iii)$,
we assume that $c_\varphi^{-1}(0)$ is not principal and that $\varphi$
satisfies
the Walters condition. According to $(ii)$ of Theorem \ref{wrpf}, there is
$h\in C(X,{\bf
R})$, $h>0$, such that $\mathcal{L}_\varphi h= \lambda h$, where
$\log\lambda=P(T,\varphi)$. Since $c_\varphi^{-1}(0)$ is not principal,
there exists $n\ge
1$ and $x\in X$ such that
$T^nx=x$ and $\varphi_n(x)=0$.
Then,
$$\lambda^n h(x)=\mathcal{L}^n_\varphi
h(x)=\sum_{T^ny=x}e^{\varphi_n(y)}h(y)> h(x). $$
Therefore $\lambda>1$ and
$P(T,\varphi)>0$. Similarly $P(T,-\beta\varphi)>0$ for all $\beta\in{\bf
R}$. According
to $(i)$, this implies that $\alpha$ has no KMS states. For
$(iv)$, assume that $\varphi$ is strictly positive. Then, the uniqueness of
$\beta$  such
that 
$P(T, -\beta\varphi) = 0$ follows by the properties of pressure
enumerated in \cite[Theorem 9.7]{w2}.  We have $P(T, 0) = h(T) > 0$
(by Corollary \ref{h(T)>0})
and $P(T, -\beta\varphi) \le h(T) - \beta\inf\varphi < 0$ for $\beta$
sufficiently
large. The continuity of $\beta\mapsto P(T, -\beta\varphi)$ and the
intermediate value
theorem give the existence of
$\beta>0$ such that
$P(T, -\beta\varphi)=0$.  The
convexity of this function implies the uniqueness of the solution.
The case for strictly negative $\varphi$ is handled similarly.
The inequality follows immediately from \cite[Theorem 9.7(ii)]{w2}.
Since the strict positivity (or negativity) of $\varphi$ implies that
$c_\varphi^{-1}(0)$ is principal, we can apply $(ii)$ to obtain the
uniqueness when
$\varphi$ also satisfies the Bowen condition.
\end{proof}

\begin{rem} If $\varphi$ is chosen to be the constant function $0$, then
$\alpha$ is the
identity automorphism group and $(iii)$ says that $C^*(X,T)$ has no tracial
state. As noted above it is well known that this C$^*$-algebra is
simple and purely infinite (see \cite{claire}).   If 
$\varphi$ is chosen to be the constant function
$1$, then
$\beta = h(T)$ the topological entropy (see \ref{main}(iv)). 
Since $\varphi$ trivially satisfies the Bowen condition, there
is a unique KMS state; this restricts to a trace on $C^*(R)$.
The associated action $\alpha$ is the usual gauge action
extended to ${\bf R}$.  This generalizes results of Olesen and Pedersen
(see \cite{op}).
\end{rem}

KMS states may fail to exist if $\varphi$ is not strictly positive (or
negative) as the next example shows.  But this is not a necessary
condition (see Example \ref{nonpos}).

\begin{ex}\label{ber2}
Let $X$ and $T$ be as in Example \ref{ber} above; recall that
$C^*(G) \cong {\mathcal O}_n$ is generated by $n$ isometries
$S_1, \dots, S_n$.   Given real numbers
$\lambda_1, \dots, \lambda_n$ there is a (unique) one-parameter
automorphism group $\alpha: {\bf R} \rightarrow {\rm Aut}({\mathcal O}_n)$
such that
$\alpha_t(S_j) = e^{it\lambda_j}S_j$ for $j = 1, 2, \dots, n$.
By \cite[Prop. 2.2]{ev} there is a KMS state for $\alpha$ at inverse
temperature $\beta$ iff
$$
1 = \sum_{j=1}^n e^{-\beta\lambda_j};
$$
in this case the KMS state is unique.  This condition holds for some
$\beta$ exactly when $\lambda_1, \dots, \lambda_n$ all have the same sign
(i.e.\ all are positive or all are negative). Let us deduce these results
from
Theorem  \ref{main}.
We define $\varphi: X \rightarrow {\bf R}$ by $\varphi(x) = \lambda_{x_0}$;
then $\varphi$
is H\"older continuous and $\alpha$ is the one-parameter automorphism group
associated to the cocycle $c_\varphi$. Note that we have
$$
\mathcal{L}_{-\beta\varphi} 1_X= (\sum_{j=1}^n e^{-\beta\lambda_j})1_X.
$$ Therefore, according to Proposition
\ref{press},
$\sum_{j=1}^n e^{-\beta\lambda_j}$ is the exponential of the pressure
$P(T,-\beta\varphi)$ and  Evans' equation is equivalent to
$P(T,-\beta\varphi)=0$. Thus,
for a function $\varphi$ of that form, a KMS state exists if and only if
$\varphi$ is strictly positive or strictly negative. Theorem  \ref{main}
applies to give
uniqueness. Note that there exist $\varphi \in C(X ; {\bf R})$ of the above
form which
do not have a constant sign but such that $c_\varphi^{-1}(0)$ is principal.
This shows
that the condition that $c_\varphi^{-1}(0)$ is principal is not sufficient
for the
existence of KMS states.
\end{ex}
\begin{ex}\label{sft2} Evans' example has been generalized to Cuntz-Krieger
algebras by J.~Zacharias in \cite{z}.   Let $A$ be an irreducible $n
\times n$ zero-one matrix $A = (A(i,j))_{1 \le i,j \le n}$ and let
$X_A$ and $T_A$ as in Example \ref{sft}. There are   
$n$ partial isometries $S_1, \dots, S_n$ generating $\,C^*(G)
\cong {\mathcal O}_A$ and one can define
$\alpha: {\bf R} \rightarrow {\rm Aut}(\mathcal{O}_A)$ just as above.
Zacharias  has shown in \cite[Prop.\ 4.3]{z} 
that there is a KMS state for
$\alpha$ if and only all the $\lambda_i$ have the same sign and in this case
the KMS
state is unique; the inverse temperature $\beta$ satisfies the
condition that the spectral radius of $D_\beta A$ is $1$, where
$D_\beta$ is the diagonal matrix with entries $e^{-\beta\lambda_i}$. This
can be deduced
from Theorem  \ref{main} just as above. We define $\varphi$ as previously.
We find that
the pressure
$P(T,-\beta\varphi)$ is the logarithm of the the spectral radius of $D_\beta
A$
and that it can take the value $0$ if and only if all the $\lambda_i$ have
the same sign.
Then, we can apply the part $(iv)$ of the theorem.
\end{ex}

\begin{rem}\label{sft3} 
  KMS states with respect to gauge automorphism groups on
Cuntz-Krieger algebras have also been studied by D.~Kerr and C.~Pinzari in
\cite[Section 7]{kp} and by R.~Exel in \cite{ex2}. Given a primitive
$n
\times n$ matrix
$A$ with zero-one elements, we define $X = X_A$ and $T = T_A$ as in Example
\ref{sft} and ${\mathcal O}_A\cong C^*(G(X_A,T_A))$. The automorphism group
$\alpha^{b,a}$, where $b\in{\bf R}$ and $a\in C(X,{\bf R})$, considered in
\cite{kp} is our
automorphism group $\alpha$ for $\varphi=b-a$. They consider the KMS problem
for
$\alpha^{b,a}$ at inverse temperature $1$. Our necessary and sufficient
condition for existence
$P(T_A,-\varphi)=0$ becomes
$b=P(T_A,a)$. In order to obtain the uniqueness of the KMS state, they
assume that the
variation of $a$ is strictly less than the entropy $P(T_A,0)$; since this
condition
implies the strict positivity of $\varphi=b-a$, their results are covered by
Theorem
\ref{main} . The existence and uniqueness of KMS states in the case where
$(X,T) = (X_A,T_A)$ and $\varphi$ is H\"older continuous and strictly positive also
appears in
\cite[Theorem 4.4]{ex2}.
\end{rem}

\begin{ex}\label{nonpos} Here is an example where the KMS state is unique
but the
potential $\varphi$ does not have a constant sign. Let $X$ and $T$ be as in
Example
\ref{ber} with
$n = 2$.  Let 
$\varphi \in C(X ; {\bf R})$ be defined by
$$
\varphi(x) =
\begin{cases}
\log(4/3) & \mbox{ if } x_0 = 0 \mbox{ and } x_1 = 1 \\
-\log 3  & \mbox{ otherwise}.
\end{cases}
$$
Then we claim that $\alpha$ as defined above has a unique KMS state which
occurs at
inverse temperature $\beta = -1$. We write $(X,T)$ as the shift associated
to the
stationary Bratteli diagram given by the matrix 
$$
\begin{pmatrix}
1 & 1 \\ 1 & 1
\end{pmatrix}.
$$ 
Then $c_\varphi$ is the stationary quasi-product cocycle
given by the matrix 
$$
\begin{pmatrix}
-\log 3 & \log(4/3)\\ -\log 3 & -\log 3
\end{pmatrix}.
$$ 
According to Proposition \ref{press} and \cite[Section
3.2]{r3}, the
pressure $P(T,-\beta\varphi)$ is the logarithm of the spectral radius of the
matrix 
$$
B_\beta =
3^\beta\begin{pmatrix}
1 & 4^{-\beta} \\ 1 & 1
\end{pmatrix}.
$$
It is a strictly increasing function from $-\infty$ to $+\infty$ and it
takes the value
$0$ at $\beta=-1$. In this example, $c_\varphi^{-1}(0)$ is
principal, $\varphi$ does not have a constant sign and $\alpha$ has a unique
KMS state.
\end{ex}

\begin{rem}
Let $T : X \rightarrow X$ be a local homeomorphism, let $\varphi \in C(X ;
{\bf R})$
and let $\alpha$ denote the action associated to the cocycle \
$c_\varphi \in Z^1(G, {\bf R})$.   Then the subalgebra $C^*(R_n)$ is
invariant
under $\alpha$ and the restriction of $\alpha$ to $C^*(R_n)$ is inner (this
is because the restriction of $c_\varphi$ to $R_n$ is a coboundary).  It
follows then that the restriction of $\alpha$ to $C^*(R)$ is
approximately inner.  It follows by \cite[Proposition 8.12.9]{p} that
$C^*(R)$ has a ground state. Moreover, one may use the methods of
\cite{r5} to show that KMS states for the restriction of $\alpha$ to
$C^*(R)$ exist at
every inverse temperature
$\beta \in {\bf R}$. This generalizes Propositions III.1.5 and III.2.9 of
\cite{r1}.
\end{rem}
\section{The crossed product}\label{cp}

We collect here some facts concerning the crossed product
$C^*(G) \rtimes_\alpha {\bf R}$.  
If $T : X \rightarrow X$ is a local
homeomorphism
and $\varphi \in C(X ; {\bf R})$, then $\tilde{T} : X \times {\bf R}
\rightarrow X \times {\bf R}$
defined by $\tilde{T}(x, t) = (x, t + \varphi(x))$ is a local
homeomorphism.  Moreover, $G(X \times {\bf R}, \tilde{T})$ may be identified
with the skew product groupoid $G(c_\varphi)$ and so
$C^*(G) \rtimes_\alpha {\bf R} \cong C^*(G(X \times {\bf R}, \tilde{T}))$
(see \cite[Example 4]{dkm}).

Note that $G$ is minimal but $G(X \times {\bf R}, \tilde{T})$ may well
not be.  We consider a necessary and sufficient condition for the minimality of
$G(X \times {\bf R}, \tilde{T})$:
For every $x \in X$, open set $U \ni x$, $t \in {\bf R}$ and $\epsilon > 0$,
there exist $y \in U$, $m, n \in {\bf N}$ such that $T^mx = T^ny$ and
$|c_\varphi(x, m - n, y) - t | < \epsilon$.
This condition is equivalent to the requirement that
$R_\infty^{\,x}(c_\varphi) = {\bf R}$ for all $x \in X$
(see \cite[Definition I.4.3]{r1}). 
Hence,  $G(X \times {\bf R}, \tilde{T})$ is minimal and thus
$C^*(G) \rtimes_\alpha {\bf R}$ is simple if and only if this condition is
satisfied (see \cite[Proposition I.4.14]{r1}). This condition may be
difficult to check but see the examples below.  But simplicity of the
crossed product also has interesting consequences.  The following
proposition is a restatement of some earlier results appropriate to
the current context:

 
\begin{prop}\label{cross}
Let $X$ be a compact metric space and let $T : X \rightarrow X$ be a
local homeomorphism which is positively  expansive and exact. 
Let $\varphi \in C(X ; {\bf R})$ and let $\alpha$ be the one-parameter
automorphism group associated to the cocycle 
$c_\varphi \in Z^1(G, {\bf R})$.  Then 
\begin{itemize}
\item[(i)] If $C^*(G) \rtimes_\alpha {\bf R}$  is simple, then it is
  stable and either infinite or projectionless.
\item[(ii)] If $\varphi$ is strictly positive (or strictly negative)
  there is a faithful lower semicontinuous trace $\tau$ on the
  crossed product $C^*(G) \rtimes_\alpha {\bf R}$  which satisfies the
  scaling property with respect to the dual automorphism, i.e.\
  $\tau\circ \hat{\alpha}_s = e^{-\beta s}\tau$ for all $s \in
  \hat{\bf R }$.  It follows that the crossed product is stably
  projectionless. 
\end{itemize}
\end{prop}
\begin{proof}
Assertion (i) follows from \cite[Props. 3, 4]{kk2} and the fact that
$C^*(G)$ is a unital simple purely infinite C*-algebra (see
\cite{claire}). 
Assertion (ii) follows from the existence of KMS states (see Theorem
\ref{main}iv), \cite[Theorem 3.2]{kk1} and \cite[Corollary 3.4]{kk1}.  
\end{proof}

\begin{rem}
In the setting of part (ii) of the above proposition it may be
worthwhile to find an explicit formula for the trace on the
crossed product constructed from a measure resulting from Proposition
\ref{press};  note that by arguing as in the proof of Theorem
\ref{main} there is a $\beta > 0$ such that $P(T, -\beta\varphi) = 0$.
Let $\mu$ be a measure on $X$ such that
$$
\mathcal{L}_{-\beta\varphi}^*\mu = \mu
$$
(see the proof of Theorem \ref{main}); note that $\mu$ has full support.
The measure $\tilde{\mu} = \mu \times e^{-\beta t}dt$ on $X \times {\bf R}$
is invariant under the groupoid $G(X \times {\bf R}, \tilde{T})$ and therefore
defines a faithful lower semicontinuous trace $\tau_\omega$ on the
crossed product: it is easy to check that it satisfies the scaling
property (see \cite[Theorem 3.2]{kk1}):
$$
\tau_\omega \circ \hat{\alpha}_s
              = e^{-\beta s}\tau_\omega.
$$ 
It follows that $C^*(G) \rtimes_\alpha {\bf R}$ is stably projectionless.
\end{rem}

\begin{exs}
In the case of Example \ref{ber2} with $\lambda_i$ all of the same sign,
Kishimoto has shown that $C^*(G) \rtimes_\alpha {\bf R}$ is simple if and
only if
the subgroup generated by the $\lambda_i$ is dense in ${\bf R}$ (see
\cite{ks}); in this case it is also stably projectionless (see
\cite[Th.\ 4.1]{kk1}).
In the case of Example \ref{sft2} with $A$ primitive and
$\lambda_i$ all of the same sign, Zacharias has shown that $C^*(G)
\rtimes_\alpha
{\bf R}$ is simple and stably projectionless if the $\lambda_i$ are linearly
independent over ${\bf Z}$ (see \cite[Prop. 4.6]{z}).

\end{exs}

\end{document}